\input amstex
\magnification 1200
\documentstyle{amsppt}
\voffset=-2truepc
\overfullrule=0pt
\vsize=9truein
\define\false{{\bold{false}}}
\define\true{{\bold{true}}}
\define\a#1{{\langle#1\rangle}}

\topmatter
\title
Seven Trees in One
\endtitle
\author
Andreas Blass
\endauthor
\address
Mathematics Dept., University of Michigan, Ann Arbor, 
MI 48109, U.S.A.
\endaddress
\email
ablass\@umich.edu
\endemail
\thanks Partially supported by NSF grant DMS-9204276.
\endthanks
\abstract
Following a remark of Lawvere, we explicitly exhibit a 
particularly elementary bijection between the set $T$ of 
finite binary trees and the set $T^7$ of seven-tuples of 
such trees.  ``Particularly elementary'' means that the 
application of the bijection to a seven-tuple of trees 
involves case distinctions only down to a fixed depth 
(namely four) in the given seven-tuple.  We clarify how 
this and similar bijections are related to the free 
commutative semiring on one generator $X$ subject to 
$X=1+X^2$.  Finally, our main theorem is that the 
existence of particularly elementary bijections can be 
deduced from the provable existence, in intuitionistic 
type theory, of any bijections at all.

\endabstract
\endtopmatter
\document

\head
Introduction
\endhead
 This paper was motivated by a remark of Lawvere 
\cite{8}, which implies that there is a particularly  
elementary coding of seven-tuples of binary trees as  
single binary trees.  In Section 1, we explicitly exhibit  
such a coding and discuss the sense in which it is   
particularly elementary.  In Section 2, we discuss the  
algebra behind this situation, which explains why   
``seven'' appears here.  Section 3 connects, in  a somewhat  
more general context, the algebraic manipulations of  
Section 2 with the elementary codings of Section 1. 
Finally, in Section 4, we prove a  meta-theorem saying 
that such particularly elementary  constructions can be 
extracted from existence proofs carried out in  the much 
more liberal context of constructive type  theory. 

Throughout this paper, we use {\sl tree\/} to  mean 
specifically a finite binary tree in which the  immediate 
successors ($=$ children) of any node are  labeled left and 
right.  Even if a node has only one child,  the child must 
have a label.  We admit the empty tree,  denoted  by 0, 
but every non-empty tree has a unique root, from  which 
every node can be reached by repeatedly passing  to 
children.  The {\sl depth\/} of a node is defined as the 
number of nodes on the path joining it to the root (so the 
root has depth 1), and the  largest depth of any node is 
the {\sl depth\/} of the tree.   (If we need to assign a 
depth to the empty tree, we assign 0.)  As this 
terminology suggests, we visualize trees as  growing 
downward, the root being at the top.  We use the  
notation $[t_1,t_2]$ for the tree consisting of a root, a left  
subtree (consisting of the left child, if any, and all its  
descendants) isomorphic to $t_1$, and a right subtree  
isomorphic to $t_2$.  For example, $[0,0]$ consists of just  
the root.  By the {\sl leftward path\/} of a tree, we mean  
the set of those nodes such that neither the node nor any ancestor 
of it is a right  child, i.e., the set consisting of the root (if 
any), its left  child (if any), its left child (if any), etc.   

\head
1. A Very Explicit Bijection
\endhead

 Our first theorem is due to Lawvere, who mentions it  
(though without proof and not quite in the same form) in  
\cite{8}. 

\proclaim{Theorem  1}  
There is a very explicit bijection between the set of all  
seven-tuples of trees and the set of all trees.
\endproclaim

 Before we prove the theorem, we must explain the  
meaning of ``very explicit,'' for without this phrase the  
theorem is trivial.  The set $T$ of all trees (as defined in 
the introduction) is clearly countably infinite  and is 
therefore in one-to-one correspondence with the  set 
$T^7$ of seven-tuples of trees (and also with $T^k$ for  
any finite $k\geq1$).  Since a bijection from $T$ to the  
set of natural numbers can be given explicitly, so can the  
bijection $T^7\to T$.  This is {\sl not\/} the content of the  
theorem.  ``Very explicit'' means, roughly, that the value  
of the bijection $f$ at a seven-tuple $\vec  
t=(t_1,\dots,t_7)$ of trees can be determined by (1)  
inspecting these seven trees down to a depth $n$ that  
depends only on $f$, not on the particular trees, (2)  
depending on these seven partial trees, constructing a  
part $p$ of $f(\vec t\/)$, and (3) taking the subtrees that  
were below depth $n$ in the $t_i$'s (and were thus  
ignored at step (1)) and attaching them to certain leaves  
of $p$.  In other words, any structure occurring below  
depth $n$ in the $t_i$'s is simply copied into $f(\vec 
t\/)$;  the ``real work'' done by $f$ involves the $t_i$'s 
only to  depth $n$. 

For a more precise definition of ``very explicit'' (which  
the impatient reader can skip for now, as it will not be  
needed until Section 3), first define a {\sl pattern\/} to be 
a  tree in which some of the leaves are labeled with 
distinct  symbols.  (This labeling has nothing to do with 
the left-right labeling that is part of the structure of 
every tree.)  An {\sl instance\/} of a pattern is obtained 
by  replacing each labeled leaf by some tree; the function  
assigning to each label the corresponding tree is called  
the {\sl substitution\/} leading to the instance.  (Note that  
the substitution can have the empty tree 0 as a value; the  
corresponding labeled node is then removed from the  
pattern when the instance is formed.)  Similarly, we  
define 7-{\sl patterns} to be seven-tuples of patterns  
with all labels distinct, and we define their instances  
under substitutions, which are seven-tuples of trees, by  
replacing every labeled node by the image of its label  
under the substitution.  A {\sl very explicit function\/} 
$f:  T^7\to T$ is given by (1) a finite indexed family 
$(\vec  p_i)_{i\in I}$ of 7-patterns such that every 
seven-tuple  $\vec t$ of trees is an instance of exactly 
one $\vec p_i$  (under exactly one substitution) and (2) a 
family  $(q_i)_{i\in I}$ of patterns, indexed by the same 
$I$, such  that each $q_i$ contains the same labels as 
$\vec p_i$.   To apply $f$ to a seven-tuple $\vec t$ of 
trees, find the  unique $\vec p_i$ and the unique 
substitution yielding  $\vec t$ as an instance, and then 
apply the same  substitution to $q_i$.   

A non-trivial example will occur in the proof of  
Theorem~1.  For a trivial but instructive example, note  
that there is a very explicit function $T^2\to T$ sending  
every pair $(t_1,t_2)$ of trees to the tree $[t_1,t_2]$  
consisting of a  root with left subtree $t_1$ and right 
subtree $t_2$. (In  the notation of the preceding 
definition, it is given by a  singleton $I$, so we can omit 
the subscripts $i$; by $\vec  p$  consisting of two 
one-point trees with the points labeled,  say, 1 and 2; and 
by $q$ a tree consisting of the root, a  left child labeled 1, 
and a right child labeled 2.)  This  very explicit map is 
one-to-one but not quite surjective,  as the empty tree is 
not in its range.  There is also a very  explicit injection in 
the opposite direction, $T\to T^2$,  sending each tree $t$ 
to the pair $(t,0)$.  One could apply  the 
Cantor-Schr\"oder-Bernstein argument to this pair of  
injections to obtain a bijection $T^2\to T$, but the result  
is not very explicit.  In fact, the bijection one obtains is  
just like the injection $T^2\to T$ described above except  
that, if the output is just a leftward path (i.e., a tree  
where no node has a right child) then the path is  
shortened by one node.  Unlike the injections, this  
bijection may need to look down to arbitrary depth in its  
input trees to determine whether the output is a leftward 
path.  We shall see in Section 3 that there is  no very  
explicit bijection $T^2\to T$. 

\demo{Proof of Theorem 1}
We shall define the desired bijection $f:T^7\to T$ by cases  
depending on the structure of its input $\vec  
t=(t_1,\dots,t_7)$ (down to depth 4).  To improve  
readability, we shall write integers $k$ in place of the  
trees $t_k$.  

 \smallskip
{\sl Case 1.}\quad
At least one of the first four trees is non-empty. Output 
the tree $[[[[[[7,6],5],4],3],2],1]$.

 {\sl Case 2.}\quad
Trees 1 through 4 are empty, but tree 5 is not, say  
$5=[5a,5b]$. Output the tree $[[[[0,7],6],5a],5b]$  

{\sl Case 3.}\quad
Trees 1 through 5 are empty, but tree 6 is not. Output the 
tree $[[[[[6,7],0],0],0],0]$ 

{\sl Case 4.}\quad
Trees 1 through 6 are empty, but the leftward path in  
tree 7 has at least 4 nodes. So $7=[[[[7a,7b],7c],7d],7e]$. 
Output the tree $[[[[[0,7a],7b],7c],7d],7e]$. 

{\sl Case 5.}\quad
Otherwise, output tree 7. \smallskip 

This $f$ clearly fits the rough description above of very  
explicit functions.  To see that it fits the precise  
description, we should list the appropriate patterns $\vec  
p_i$ and $q_i$.   There are eleven $\vec p$'s and eleven  
corresponding $q$'s --- four for Case 1 (according to  
which of trees 1 through 4 is the first nonempty one), one  
each for Cases 2, 3, and 4, and four for Case 5 (according  
to whether the leftward path of tree 7 has 0, 1, 2, or 3  
nodes).  We refrain from explicitly exhibiting the results 
of this splitting of cases.

We must still verify that we have defined a bijection, i.e.,  
that every tree arises exactly once as the output of this  
construction.  For this purpose, note that the leftward  
path of the output has length 

$\geq 6$ in Case 1, 

4 in Case 2, 

$\geq 6$ in Case 3 (remember that tree 6 is non-empty  
here), 

5 in Case 4, and 

$\leq 3$ in Case 5. 

This shows that no two cases, with the possible exception  
of Cases 1 and 3, could produce the same tree as output.   
The possible exception is not a real exception, because in  
Case 1 at least one of the first four nodes on the leftward  
path has a right successor (because those successor  
subtrees are 1 through 4, of which at least one is  
non-empty by case hypothesis) whereas in Case 3 this is  
not the case.  Thus, given an arbitrary tree $t$, we can  
determine the unique case that might produce it as  
output.  Once this is done, it is straightforward to check  
that each case actually does produce all the trees thereby  
assigned to it, and that it produces them exactly once  
each.
\qed\enddemo 

Theorem~1 would remain true if we extended the concept 
of tree to allow infinite trees, and the proof would be 
unchanged.  The reason is that a very explicit bijection 
can be applied to infinite trees just as to finite trees, since 
below a certain depth it merely copies the input into 
appropriate places in the output.  More unusual notions 
of ``tree'' can be handled similarly.  For example, we 
could allow trees to have finitely many infinite paths; we 
could allow infinite paths only if, beyond some node on a 
path, all further nodes are left children; we could fix a set 
$C$ of ``colors'' and allow infinite trees in which every 
infinite path is assigned a color from $C$ (and we could 
impose continuity requirements on the coloring); etc.  In 
each case, Theorem~1 remains true as long both 
occurrences of ``trees'' are interpreted the same way. 

\head
2. Algebra
\endhead
The proof of Theorem 1 raises at least two questions.   
Where did it come from? And why seven?  The reader is  
invited to try some ``easier'' numbers in place of seven,  
say five or two; no analogous proof will be forthcoming,  
so there is indeed something special about seven.  (Of  
course there {\sl is\/} an analogue for thirteen.  Given  
thirteen trees, apply Theorem 1 to code the first seven as  
a single tree, and then code this tree with the other six of  
the original inputs.  The same trivial observation handles  
any number congruent to 1 modulo 6.)   

To see why seven is special, we first give an argument to 
establish Theorem 1 in  the style of eighteenth-century 
analysis, where  meaningless computations (e.g., 
manipulating divergent  series as though they converged 
absolutely and  uniformly) somehow gave correct results.  
This argument  begins with the observation that a tree 
either is 0 or  splits naturally into two subtrees (by 
removing the root).   Thus the set $T$ of trees satisfies 
$T=1+T^2$.  (Of course  equality here actually means an 
obvious isomorphism.  Note that the same equation holds 
also for the variant notions of tree mentioned at the end 
of Section~1.)   Solving this quadratic equation for $T$, we 
find $T=\frac12\pm i\frac{\sqrt3}2$.  (The reader who 
objects  that this is nonsense has not truly entered into 
the  eighteenth-century spirit.)  These complex numbers 
are  primitive sixth roots of unity, so we have $T^6=1$ 
and  $T^7=T$.  And this is why seven-tuples of trees can 
be coded as single trees. 

Although this computation is nonsense, it has at least the  
psychological effect of suggesting that something like  
Theorem 1 has a better chance of being true for seven  
than for five or two.  To improve the effect from  
psychological to mathematical, we attempt to remove the  
nonsense while keeping the essence of the computation.  
Incidentally, Lawvere's remark that led to this paper was  
phrased in terms of such a meaningless computation  
giving a correct result: ``I was surprised to note that an  
isomorphism $x=1+x^2$ (leading to complex numbers as  
Euler characteristics if they don't collapse) always  
induces an isomorphism $x^7=x$'' \cite{8} p. 11.  This remark was an
application of the method of Schanuel \cite{10} for developing
objective meaning for such computations.

One might hope for a meta-theorem to the effect that,  
when such a meaningless computation leads from a  
meaningful equation (like $T=1+T^2$) to another  
meaningful equation (like $T^7=T$) then the latter  
honestly follows from the former.  This would be  
somewhat analogous to Hilbert's program for making  
constructive  sense of infinitary mathematics by showing 
that, when a  detour through the infinite leads from one 
finitistically  meaningful statement to another, then the 
latter honestly  (finitistically) follows from the former.  
Unfortunately,  our situation (like Hilbert's) is more 
subtle.  After all,  $T^6=1$  is a meaningful equation, 
saying that there is  only one six-tuple of trees, which is 
false.  So any  rehabilitation of the computation above 
must, in  particular, explain why the ultimate conclusion 
$T^7=T$ is  right while the penultimate $T^6=1$ is not. 

As a first step toward rehabilitation, we observe that we  
can avoid complex numbers by working entirely within  
the ring of polynomials with integer coefficients.  (The  
possibility of such a step is guaranteed by general facts  
about polynomial rings and ideals, but here is the step  
explicitly.)  We simplify $T^7$ by repeatedly reducing its  
degree, using the equation $T=1+T^2$ in the form $T^2=T- 
1$.  The result is  $$
T^7=T^6-T^5=-T^4=-T^3+T^2=T.
$$
The complex numbers are gone, but the computation is  
still meaningless (when we remember that $T$ is a set)  
because of negative coefficients, and we could still get  
$T^6=1$ just by reducing all exponents by one. 

The next step is to eliminate the negative terms by  
adding new terms ($T^5$, $T^4$, and $T^3$) to cancel  
them.  The resulting calculation, which takes place in the  
semiring $\Bbb N[T]/(T=1+T^2)$, reads $$\align
T^7+T^5+T^4+T^3&=T^6+T^4+T^3=T^5+T^3\\
&=T^5+T^4+T^2=T^5+T^4+T^3+T.
\endalign$$
Everything here is meaningful and correct when $T$ is  
interpreted as  the set of trees (rather than a formal 
indeterminate  subject to $T=1+T^2$) and equality is 
interpreted as  obvious isomorphism (or, better, as very 
explicit  bijection).  Unfortunately, instead of getting a 
very  explicit bijection from $T^7$ to $T$, we have the 
extra  terms $T^5+T^4+T^3$ on both sides.  Can we get rid 
of  them? 

There is a general way to convert a bijection $f:A+X\to  
B+X$ into a bijection $f':A\to B$ when $X$ is finite,  
namely the Garsia-Milne involution principle \cite{6}.   
The idea is to define $f'(a)$ by first applying $f$ to $a$; if  
the result is in $B$, accept it as the output for $f'$; if not,  
then it is in $X$, so we can apply $f$ again;  if the result 
is in $B$, accept it as the output for $f'$; if  not, then it is 
in $X$, so we can apply $f$ again; and so  forth until we 
finally get an output in $B$.  The finiteness  of $X$ is used 
to ensure that the process terminates.  In  our situation, 
$X=T^5+T^4+T^3$ is not finite.  One can try to  apply  the 
Garsia-Milne construction anyway, hoping that the  
process terminates even though no finiteness forces it to.   
Unfortunately, it does not terminate. 

We can do better by considering the last computation  
displayed above.  It is in two parts.  In the first part,  
$T^7+T^5+T^4+T^3$ is simplified to $T^5+T^3$; that is,  
$T^7+T^4$ is eliminated.  In this elimination, the presence  
of $T^5$ was essential as a sort of catalyst, since the $T^5$  
is removed at the first step and restored at the second.   
$T^3$, on the other hand is irrelevant to this half of the  
calculation.  Clearly, the same argument shows that, in  
any polynomial in $\Bbb N[T]/(T=1+T^2)$, we can delete  
or introduce $T^{k+3}+T^k$ provided $T^{k+1}$ is also  
present as a catalyst.  The second half of the displayed  
computation introduces $T^4+T$, using $T^3$ as a catalyst;  
the same argument would allow us to introduce or delete  
$T^{k+3}+T^k$ provided $T^{k+2}$ is also present as a  
catalyst.  Thus, we can delete or introduce two powers of  
$T$ whose exponents differ by exactly three, provided  
one of the two intervening powers is also present as a  
catalyst. 

This result can be significantly improved by noticing that  
{\sl any\/} positive power of $T$ can serve as a catalyst;  
it need not be between the two that are being deleted or  
introduced.  To see this, suppose we want to delete or  
introduce $T^{k+3}+T^k$, and we have another positive  
power $T^r$ present in the polynomial.  If $r$ is $k+1$ or  
$k+2$, then we already know $T^r$ can serve as the  
desired  catalyst.  If not, then use $T=1+T^2$ to replace 
$T^r$ with  $T^{r-1}+T^{r+1}$.  So now we have two 
powers of $T$, in  one of which the exponent is closer to 
the desired $k+1$  or $k+2$.  Apply the same procedure 
to that power,  leaving the other one alone, and  repeat 
the process until finally you get the desired  catalyst.  Use 
it to delete or introduce $T^{k+3}+T^k$, and  then reverse 
the previous steps to recover the original $T^r$ by 
reassembling the terms into which the first part of the 
procedure decomposed it. (Note that we needed $r$ to be 
strictly positive in order  to  replace $T^r$ with $T^{r-
1}+T^{r+1}$.)  

Now we can, in $\Bbb N[T]/(T=1+T^2)$, convert $T^7$ into  
$T^7+T^4+T$ (introducing $T^4+T$ with $T^7$ as catalyst)  
and then into $T$ (deleting $T^7+T^4$ with $T$ as  
catalyst). We refrain from writing out explicitly the 
computation in  $\Bbb N[T]/(T=1+T^2)$ obtained by the 
method of the  preceding paragraph.  It consists of 
twenty steps: four to  convert $T^7$ into a usable catalyst 
$T^3$ (plus extra  terms), two to introduce $T^4+T$ using 
this catalyst, four  to collect the catalyst and extra terms 
back into a $T^7$,  four to convert $T$ into a usable 
catalyst $T^5$ (plus  extra terms), two to delete $T^7+T^4$ 
using this catalyst,  and four to collect the catalyst and 
extra terms back into  a $T$.  Notice that the proof of 
$T^7=T$ does not become a  proof of $T^6=1$ if we reduce 
the exponents; unlike $T$,  1 cannot serve as a catalyst. 

We can now answer the first of the questions at the  
beginning of this section --- where did the bijection in  
the proof of Theorem~1 come from?  It came from this  
twenty-step computation.  Wherever the computation  
used $T=1+T^2$, apply the obvious bijection between the  
sets $T$ and $\{0\}\cup T^2$, and the whole computation  
will give a composite bijection which is the one used to  
prove Theorem~1 (except that I interchanged left and  
right once, to make the description of the proof easier). 

As for the second question --- why seven? --- we have a  
partial answer.  The technique that led to the proof of  
Theorem~1 works only for numbers $k$ that are  
congruent to 1 modulo 6, because it is only for these  
numbers that $T^k=T$ is true in the semiring $\Bbb  
N[T]/(T=1+T^2)$, or even in the rings $\Bbb  
Z[T]/(T=1+T^2)$ or $\Bbb C[T]/(T=1+T^2)$, i.e., because it  
is only for these values of $k$ that $T^k=T$ is satisfied by  
the  roots $\frac12\pm i\frac{\sqrt3}2$ of $T=1+T^2$.  A 
more  complete answer to the second question would 
involve  showing that very explicit bijections between 
$T^k$ and  $T$ can  exist only when the corresponding 
polynomials are equal  in $\Bbb N[T]/(T=1+T^2)$.  This 
will be done, in greater  generality, in the next section. 

Because of the importance of the semiring $\Bbb 
N[T]/(T=1+T^2)$ in computations like the preceding ones, 
we present a normal form for its elements.

 \proclaim{Theorem  2}
Every element of $\Bbb N[T]/(T=1+T^2)$ is uniquely  
expressible in the form $a+bT^2+cT^4$ with non-negative  
integer coefficients $a$, $b$, and $c$ such that either at  
least one coefficient is 0 or else $a\geq b=c=1$. 
\endproclaim

 \demo{Proof}
Every element of the semiring $\Bbb N[T]/(T=1+T^2)$ is a  
polynomial in $T$ with non-negative integer coefficients.  
We show how to simplify such a polynomial to the form  
claimed in the theorem. Since $T^7=T$, all terms of degree 
7 or more can be  reduced to lower degree, so we may 
assume the  polynomial has degree at most 6.  The degree 
can be  reduced to at most 4 by means of the equations 
$$
T^6=T^5+T^7=T^5+T \text{\qquad and\qquad}
T^5=1+T^3+T^5=1+T^4,
$$
where $1+T^3$ was introduced in the second equation 
with catalyst $T^5$.
Furthermore, we can eliminate all terms of odd degree,  
since $T=1+T^2$ and $T^3=T^2+T^4$.  So our polynomial  
has the form $a+bT^2+cT^4$.  We check next that the  
coefficient restrictions in the theorem can be enforced. 

For this purpose, consider the polynomial $Q=1+T^2+T^4$.  In the
presence of a catalyst (any positive power of $T$), $Q$ vanishes, for
it equals $T+T^4$ which a catalyst can delete.  Now consider our
general polynomial $a+bT^2+cT^4$. If at least one of the coefficients
is 0, then it has the desired form, so suppose all three coefficients
are at least 1.  Then the polynomial is $Q+(a-1)+(b- 1)T^2+(c-1)T^4$.
If either $b-1$ or $c-1$ is positive, we have a catalyst and can
remove $Q$.  Repeating the process of splitting off a $Q$ and deleting
it as long as a catalyst is available, we ultimately get either a
polynomial of the form $a+bT^2+cT^4$ with a zero coefficient (if, at
some stage, we cannot split off a $Q$) or one of the form $d+Q$ with
$d\in\Bbb N$ (if, at some stage, we have split off a $Q$ but have no
catalyst to remove it).  This completes the proof that every
polynomial can be put in the required form; it remains to show
uniqueness.

If two expressions of the form described in the theorem,  
say $a+bT^2+cT^4$ and $a'+b'T^2+c'T^4$, are equal in  
$\Bbb N[T]/(T=1+T^2)$, then they are also equal as  
complex numbers when $T$ is interpreted as the  
primitive sixth root of unity $t=\frac12+i\frac{\sqrt3}2$  
(since $\Bbb C$ is a commutative semiring, and $t=1+t^2$, 
and $\Bbb N[T]/(T=1+T^2)$ with $T$ is the initial  
commutative semiring with a solution of $T=1+T^2$).   
$t^2$ is a primitive cube root of unity; its minimal  
polynomial is $1+z+z^2$.  So this minimal polynomial  
would have to divide $(a-a')+(b-b')z+(c-c')z^2$, which  
means that $a-a'$, $b-b'$, and $c-c'$ are all equal.  If they  
are all zero, then $a+bT^2+cT^4$ and $a'+b'T^2+c'T^4$ are  
identical as expressions, which is what we needed to  
prove.  So suppose that the differences $a-a'$, $b-b'$, and  
$c-c'$ are all positive.  (If they are all negative,  
interchange primed and unprimed in what follows.)  A  
fortiori, $a$, $b$, and $c$ are all positive, so  
$a+bT^2+cT^4$, being of the form in the theorem, has  
$b=c=1$ and $a\geq1$, i.e., it is  $(a-1)+Q$.  As $a-a'$, $b-
b'$ and $c-c'$ are equal and  positive, we must have 
$b'=c'=0$ and $a'=a-1$. 

Our task is thus reduced to showing that $a'$ and $a'+Q$ 
are not equal in  $\Bbb N[T]/(T=1+T^2)$.  For this 
purpose, notice that the  set consisting of the natural 
numbers and the cardinal  $\aleph_0$ is a commutative 
semiring (with the usual  addition and multiplication of 
cardinal numbers) in which  $\aleph_0=1+{\aleph_0}^2$.  
So an equation $a'=a'+Q$ in  $\Bbb N[T]/(T=1+T^2)$ would 
imply the same equation in  this cardinal semiring with 
$T$ interpreted as  $\aleph_0$.  But in this interpretation 
such an equation  is false, since $a'$ is interpreted as the 
integer $a'$ while  $a'+Q$ is interpreted as $\aleph_0$.  
So the equation  cannot hold in $\Bbb N[T]/(T=1+T^2)$. 
\qed\enddemo  

\head
3. Algebraic Equivalence and Very Explicit Bijections
\endhead

 Let us call two polynomials, $P(X)$ and $Q(X)$, with  
non-negative integer coefficients {\sl algebraically  
equivalent\/} if the equation $P(X)=Q(X)$ holds in the  
semiring $\Bbb N[X]/(X=1+X^2)$.  We call the 
indeterminate $X$ rather than $T$ to avoid confusion, 
since we  shall need to discuss it and the set $T$ of trees 
in the  same context.   

A polynomial $P(X)\in\Bbb N[X]$ has a natural  
interpretation as an operation $S\mapsto P(S)$ on sets  
(well-defined up to canonical bijections); sums and  
products of polynomials correspond to disjoint unions and  
Cartesian products of sets.  (We can interpret the  
numerical coefficients in $P(X)$ as canonically chosen sets  
of the corresponding cardinalities, or we can eliminate  
these coefficients in favor of sums of repeated terms.) 
Thanks to the canonical bijection from $T$ to $1+T^2$,  
algebraically equivalent polynomials, when applied to the  
set $T$ of trees, give canonically isomorphic sets.  In fact,  
we shall show, as half of the next theorem, that the  
canonical isomorphism is very explicit in the sense  of 
Section 1. 

Let us call two polynomials $P(X)$ and $Q(X)$ {\sl  
combinatorially equivalent\/} if there is a very explicit  
bijection from $P(T)$ to $Q(T)$. Of course, we must 
extend the definition of ``very explicit,'' given in 
Section~1 for the special case of $T^7$ and $T$, to the 
general case of $P(T)$ and $Q(T)$, but this is 
straightforward.  If $P(X)=\sum_kc_kX^k$, then we 
regard $P(T)$ as the set of tagged $k$-tuples of trees, 
$\vec t=(t_1,\dots,t_k,\tau)$ where $k$ ranges over the 
same finite set as in the sum defining $P(X)$ and where 
the tag $\tau$ is an integer in the range $1\leq\tau\leq 
c_k$.  For brevity, we call such a tagged $k$-tuple a 
$P$-tuple.  A $P$-{\sl pattern\/} is a similarly tagged 
$k$-tuple of patterns (in the sense defined in Section~1) 
in which all the labels are distinct, and the notion of an 
{\sl instance} of a pattern with respect to a substitution is 
defined just as in Section~1.  A {\sl very explicit 
function\/} $f$ from $P(T)$ to $Q(T)$ is given by an 
indexed family $(\vec p_i)_{i\in I}$ of $P$-patterns and a 
similarly indexed family $(\vec q_i)_{i\in I}$ of 
$Q$-patterns such that, for each $i\in I$, the same labels 
occur in $\vec p_i$ as in $\vec q_i$ and such that every 
element $\vec t$ of $P(T)$ is an instance of a unique 
$\vec p_i$ (under a unique substitution).  Then $f(\vec 
t\/)$ is defined as the result of applying the same 
substitution to $\vec q_i$ (for the same $i$). 

\remark{Remark}
If the very explicit function $f$ is a bijection from $P(T)$ 
to $Q(T)$, then every element of $Q(T)$ occurs exactly 
once as an instance of a $\vec q_i$ (in the notation of the 
preceding definition).  It follows that the inverse bijection 
is also very explicit, as we can interchange the roles of 
the $\vec p_i$ and the $\vec q_i$.  Thus, it makes no 
difference whether ``very explicit bijection'' is 
interpreted in the obvious way as ``very explicit function 
that happens to be a bijection'' or as ``very explicit 
function with very explicit two-sided inverse.''

Our interest will be focused on very explicit bijections, 
and for the study of these our definition of ``very 
explicit'' seems adequate.  If we were to study very 
explicit functions in general, it would be advisable to 
liberalize the definition of ``very explicit'' by allowing 
labels to be repeated in a $\vec q_i$ and allowing a label 
to occur in $\vec p_i$ without occurring in $\vec q_i$.  In 
either of these cases, the very explicit function cannot be 
a bijection.
\endremark

\proclaim{Theorem  3}
Two polynomials $P(X), Q(X)\in\Bbb N[X]$ are  
combinatorially equivalent if and only if they are  
algebraically equivalent.
\endproclaim 

\demo{Proof}
Suppose first that $P(X)$ and $Q(X)$ are algebraically  
equivalent.  Since $\Bbb N[X]$ is the free commutative  
semiring on one  generator, its quotient $\Bbb 
N[X]/(X=1+X^2)$ is the initial  algebra in the variety $\Cal 
V$ of commutative semirings  with a  distinguished 
element $X$ subject to $X=1+X^2$.   Algebraic equivalence 
therefore means that the equation  $P(X)=Q(X)$ is 
satisfied by the whole variety $\Cal V$.   (Note that here 
$X$ is a constant for the distinguished  element of a $\Cal 
V$-algebra.)  By Birkhoff's theorem  \cite{3, 
Theorem~14.19}, there is an equational deduction of 
$P(X)=Q(X)$  from $X=1+X^2$ and the axioms for 
commutative  semirings, using as  rules of inference only 
substitution of equals for equals  and substitution of 
terms for variables.  All the  substitutions for variables 
can be done at the beginning,  as substitutions into 
axioms.  So we get a deduction of  $P(X)=Q(X)$ from 
$X=1+X^2$ and variable-free instances  of the 
commutative semiring axioms, using only  substitution of 
equals for equals. 

Consider the following property of equations: When the  
two sides are applied as operations to $T$, the two  
resulting sets have a bijection between them such that  
both the bijection and its inverse are very explicit.  It is 
trivial to check that the equation  $X=1+X^2$ and all 
variable-free instances of the axioms of  commutative 
semirings enjoy this property and that the  property is  
preserved by substitution of equals for equals.  It  
therefore follows that $P(X)=Q(X)$ has this property,  
which implies combinatorial equivalence. 

For the converse, suppose that $P(X)$ and $Q(X)$ are  
combinatorially equivalent.  Specifically, let $f$ be a very 
explicit bijection from $P(T)$ to $Q(T)$, and let $(\vec 
p_i)_{i\in I}$ and $(\vec q_i)_{i\in I}$ be as in the 
definition of ``very explicit'' preceding the theorem.

It will be convenient, both for this proof and for 
Section~4, to introduce certain standard collections of 
patterns.  (We temporarily deal with patterns in the 
sense of Section~1, single trees; we will return to tuples 
and $P$-patterns later.)  For each positive integer $n$, let 
$S_n$ be the set of patterns $p$ of depth $\leq n+1$ such 
that all nodes at level $n$ have exactly two children,
all nodes at level $n+1$ are labeled (with distinct labels, 
as required by the definition of pattern), and no nodes at levels
$\leq n$ are labeled. If $p\in S_n$ 
then its instances are exactly those trees that are 
identical with $p$ down to depth $n$, with no constraints 
on what (if anything) happens at greater depth. In 
keeping with this, we can define $S_0$ to consist of just 
one pattern, which consists of a single, labeled node.  
Clearly, for every $n$, every tree is an instance of a 
unique member of $S_n$.  

If $r$ is any pattern of depth $< n$, then we can associate 
to it a set $\hat r\subseteq S_n$ having collectively the 
same instances as $r$.  One can produce $\hat r$ from $r$ 
by the following step-by-step procedure, which we call 
{\sl developing\/} $r$ (to depth $n$).  Choose any labeled 
leaf in $r$ and replace $r$ with  two new patterns, $r'$ in 
which this leaf has been deleted, and $r''$ in which this 
leaf has been given two labeled children (and has lost its 
own label, as it is no longer a leaf).  Apply the same 
construction to $r'$ and $r''$, using labeled leaves at 
depths $\leq n$.  Iterate the process until all the patterns 
have their labels at level $n+1$.  At this stage, they are 
easily seen to be in $S_n$.  Furthermore, the instances of 
$r$ are precisely the instances of $r'$ and those of $r''$ 
(according to whether the substitution gives the label of 
the chosen leaf the value 0 or not).  And no tree is an 
instance of both $r'$ and $r''$.  It follows inductively that, 
at each stage of the development, the patterns have 
disjoint sets of instances and the union of these sets 
contains precisely the instances of $r$.  Thus, the final set 
of patterns obtained by this development is $\hat r$.

We need to extend the preceding concepts from patterns 
to $P$-patterns (and $Q$-patterns).  By $S_n(P)$ we 
mean the set of all $P$-patterns whose component 
patterns lie in $S_n$.  (Recall that a $P$-pattern is a 
tagged $k$-tuple of patterns, so this makes sense.)  Every 
$P$-tuple of trees is an instance of a unique pattern in 
$S_n(P)$ (for each fixed $n$).  $P$-patterns can be 
developed componentwise; that is, the development of 
$(r_1,\dots,r_k,\tau)$ to depth $n$ ($>$ the depths of all 
the $r_j$'s) consists of the patterns 
$(z_1,\dots,z_k,\tau)\in S_n(P)$ with each $z_j$ in $\hat 
r_j\subseteq S_n$.  Again, the members of the 
development of $\vec r$ have disjoint sets of instances 
and the union of these sets contains precisely the 
instances of $\vec r$.

Let us return to the families $(\vec p_i)_{i\in I}$ and 
$(\vec q_i)_{i\in I}$ determining our very explicit 
bijection $f$.  Fix $n$ greater than the depths of all the 
patterns in all these $P$- and $Q$-patterns.  Since every 
$P$-tuple is an instance of exactly one $\vec p_i$, it 
follows from the preceding discussion that the sets $\hat 
p_i\subseteq S_n(P)$ obtained by developing $\vec p_i$ 
to depth $n$ are disjoint and their union is $S_n(P)$.

By the {\sl weight\/} of a pattern $r$, we mean the monomial
$X^l\in\Bbb N[X]/(X=1+X^2)$, where $l$ is the number of labels in $r$.
By the {\sl weight\/} of a (tagged) tuple of patterns (e.g., a $P$- or
$Q$-pattern), we mean the product of the weights of its component
patterns, i.e., $X$ with exponent the total number of labels in the
tuple.  By the {\sl weight\/} of a set of patterns or of (tagged)
tuples of patterns, we mean the sum of the weights of its members.

One step in the development of a pattern $r$ replaces it 
by a set of two patterns $r'$ and $r''$ where, if $r$ has 
weight $X^l$, then $r'$ has weight $X^{l-1}$ (as one 
labeled leaf was deleted) and $r''$ has weight $X^{l+1}$ 
(as one labeled leaf lost its label but was given two 
labeled children).  So the weight of the set obtained is 
$X^{l-1}+X^{l+1}=X^l$, since the weights are in  a semiring 
where $1+X^2=X$.  So one step of development leaves the 
weight unchanged.  It follows by induction that 
development of a pattern to depth $n$ leaves the weight 
unchanged; $r$ and $\hat r$ have the same weight.  It 
further follows easily that the same applies to (tagged) 
tuples of patterns.  

These observations imply that the set $S_n$ has weight 
$X$, because $S_0$ trivially has weight $X$ and $S_{n+1}$ 
is obtainable by developing $S_n$.  It follows that 
$S_n(P)$ has weight $P(X)$, because the contribution to 
the weight from the $k$-tuples with a particular tag 
$\tau$ is 
$$\align
&\sum_{\vec r\in (S_n)^k} \text {weight}(\vec r\/)=
\sum_{\vec r\in (S_n)^k} \prod_{j=1}^k \text 
{weight}(r_j)=\\
&\prod_{j=1}^k\sum_{r\in S_n} \text {weight}(r)=
\prod_{j=1}^k\text {weight}(S_n)=
X^k.
\endalign$$
Since, as we saw above, $S_n(P)$ is obtainable by 
developing the set of $\vec p_i\/$'s, this set also has 
weight $P(X)$.  Similarly, the set of $\vec q_i\/$'s has 
weight $Q(X)$.  But each $\vec p_i$ contains exactly the 
same labels as the corresponding $\vec q_i$, so these two 
sets have the same weight.  Therefore, $P(X)=Q(X)$ in 
$\Bbb N[X]/(X=1+X^2)$.
\qed\enddemo 

It follows immediately from Theorem~3 that, if $P(X)$ and $Q(X)$ are
combinatorially equivalent, then $P(t)=Q(t)$ where
$t=\frac12+i\frac{\sqrt3}2$.  Thus, the complex number $P(t)$ serves
as an invariant of the combinatorial equivalence class of $P(X)$.
This is what Lawvere referred to as ``complex Euler characteristics''
in the passage, quoted in Section~2, that motivated this paper.  In
fact, Theorem~3 and the proof of Theorem~2 show that the semiring of
combinatorial equivalence classes of polynomials is embedded in the
product of the complex field and the semiring of cardinals
$\leq\aleph_0$ (by sending the indeterminate $X$ to $\a{t,\aleph_0}$).
Lawvere has pointed out that for this purpose one could replace the
cardinal semiring with a three-element semiring, for all the finite,
non-zero cardinals can be identified without damaging the embedding.
The resulting three elements form the system of ``dimensions''
associated to our problem by Schanuel's general construction
\cite{10}, so that the present solution of the word problem shows in
particular that again in our case ``Euler characteristic and
dimension'' are jointly injective.

\head
4. Constructive Set Theory
\endhead

We show in this section that the algebraic equivalence of 
$P(X)$ and $Q(X)$ can be deduced from the mere 
provability of ``there is a bijection from $P(T)$ to $Q(T)$'' 
provided this provability is from sufficiently restricted 
assumptions.  The restrictions we need are two: The 
underlying logic is constructive and the only assumption 
about $T$ is the existence of a bijection $T\to1+T^2$.  
(Actually, the second restriction can be relaxed a bit by 
allowing a stronger assumption about $T$.) On the other 
hand, we allow the use of higher-order logic, so many 
set-theoretic methods are available.  

More precisely, let $\Cal L$ be a higher-order theory in 
the sense of \cite{4} or a local set theory in the sense of 
\cite{1}, generated by a natural number object and an 
additional ground type $T$ subject to the axiom ``there is 
a bijection $1+T^2\to T$.''  Alternatively, we could let 
$\Cal L$ be intuitionistic Zermelo-Fraenkel set theory 
augmented with a constant $T$ and the same axiom.

\proclaim{Theorem  4}
Let $P(X)$ and $Q(X)$ be polynomials with non-negative 
integer coefficients.  Suppose it is provable in $\Cal L$ 
that there is a bijection $P(T)\to Q(T)$.  Then $P(X)$ and 
$Q(X)$ are algebraically equivalent.
\endproclaim

Before proving the theorem, we make several remarks.
First, the converse of the theorem is easy to prove.  By 
Theorem~3, algebraic equivalence implies the existence of 
a very explicit bijection $P(T)\to Q(T)$ when $T$ is the 
set of trees.  But the very-explicitness makes it possible 
to apply the bijection to arbitrary sets for which a 
bijection $1+T^2\to T$ is given, and this application can 
be carried out in constructive set or type theory.  
Alternatively, we can proceed as in the proof of 
Theorem~3, considering for all equations $P(X)=Q(X)$ the 
property ``it is provable in $\Cal L$ that there is a 
bijection $P(T)\to Q(T)$,'' noticing that this property is 
enjoyed by the equation $1+X^2=X$ and by all 
variable-free instances of the axioms for commutative 
semirings, noticing further that the property is preserved 
by substitution of equals for equals, and concluding that 
the property holds of all equations true in $\Bbb 
N[X]/(X=1+X^2)$.

The remaining remarks are intended to justify the 
restrictions we place on the logic and the assumptions on 
$T$ in the theory $\Cal L$.  

If we allowed full classical set theory, with the axiom of 
choice, then the assumption $T\cong1+T^2$ implies that 
$T$ is infinite and therefore $T$, $1+T$, $T+T$, and $T^2$ 
all have the same cardinality.  It follows that every 
non-constant polynomial is equivalent, in the sense of 
provable bijection, to $T$.  In other words, with this 
stronger set theory, the corresponding notion of algebraic 
equivalence would be equality not in $\Bbb 
N[X]/(X=1+X^2)$ but in $\Bbb N[X]/(X=1+X=X+X=X^2)$, a 
semiring isomorphic to the $\Bbb N\cup\{\aleph_0\}$ 
example used at the end of the proof of Theorem~2.

If we work in classical set theory without the axiom of 
choice, so that addition and multiplication of infinite 
cardinals are no longer trivial, we still get the same 
conclusion with a bit more work.  From $T\cong1+T^2$ 
and its immediate consequence $T^2\cong T+T^3$, we 
infer that each of $T$ and $T^2$ can be embedded in the 
other.  By the Cantor-Schr\"oder-Bernstein Theorem, 
whose proof does not require the axiom of choice, we 
have $T\cong T^2$.  Then from $1\leq2\leq T$ (where 
$\leq$ means embeddability, the usual inequality 
relation on cardinals) we get $T\leq T+T\leq T^2\leq T$, 
so another application of the Cantor-Schr\"oder-Bernstein 
Theorem gives $T\cong T+T$.  

If we use intuitionistic rather than classical logic, then the 
Cantor-Schr\"oder-Bernstein Theorem is no longer 
available and, as Theorem~4 shows, the argument in the 
preceding paragraph breaks down.  Even in intuitionistic 
logic, however, if we assume that $T$ is the set of trees 
(rather than some arbitrary set with a bijection $1+T^2\to 
T$) then the argument in the preceding paragraph works, 
since the bijections produced by the 
Cantor-Schr\"oder-Bernstein Theorem can, in this case, be 
constructively defined.  For example, we described, just 
before the proof of Theorem~1, a bijection $T^2\to T$, and 
that description is intuitionistically legitimate.  The main 
point here is that the case distinction, whether a tree is 
just a leftward path, is decidable because the tree is 
finite.  There is a similarly constructive bijection $T+T\to 
T$ when $T$ is the set of finite trees.  It sends any tree 
$t$ from the first copy of $T$ to $[0,t]$, and it sends any 
$t$ from the second copy of $T$ to $[t,0]$ unless $t$ is of 
the form 0 or $[p,q]$ or $[[p,q],0]$ or $[[[p,q],0],0]$ or 
\dots with $p$ and $q$ both $\neq0$, in which case it 
sends $t$ to $t$.
Again, it is the finiteness of the trees that makes the case 
distinction decidable and the definition constructively 
correct.

\demo{Proof of Theorem 4}
We begin by giving a more useful description of sets $T$ 
with bijections $f:1+T^2\to T$.  The bijection is 
determined by specifying a distinguished element of $T$ 
and a binary operation on $T$; the distinguished element 
is the value of $f$ at the unique element of 1, and the 
operation is the restriction of $f$ to $T^2$.  To match the 
notation used earlier for trees, we write the distinguished 
element as 0 and the operation as $[-,-]$.  Thus, $T$ (with 
the structure $f$) is an algebra with one constant and one 
binary operation.  That $f$ is a bijection means that this 
algebra must satisfy the following system $\Cal T$ of 
axioms, which, for later convenience, we write as 
geometric sequents (as defined in \cite {1, page 250} or 
\cite{7, Section~6.5}), indeed, whenever possible, as 
universal Horn formulas.
\roster
\item $0=[x,y]\implies\false$
\item $[x,y]=[x',y']\implies x=x'$
\item $[x,y]=[x',y']\implies y=y'$
\item $\true\implies x=0\lor\exists y\exists z\,x=[y,z]$
\endroster
The algebra of finite trees is initial in the variety of 
algebras of signature $\{0, [-,-]\}$, and, since it satisfies 
the axioms of $\Cal T$, it is also initial in the category of 
models of $\Cal T$.  Any model $M$ of $\Cal T$ can be 
regarded as an algebra of generalized trees, in that it has 
an element 0 corresponding to the empty tree, and all its 
other elements are uniquely of the form $[y,z]$ and can 
therefore be pictured as consisting of a root with two 
subtrees, $y$ and $z$, attached to it.  Among such 
algebras are, for example, the collections of trees in any 
of the generalized senses mentioned at the end of 
Section~1.

We shall prove the theorem by constructing a specific 
topos model of $\Cal L$.  The hypothesis of the theorem 
says that in this model there must be a bijection $P(T)\to 
Q(T)$, and an analysis of what this means will lead to the 
desired conclusion.  Perhaps the most natural topos 
model of $\Cal L$ is the classifying topos (\cite {7, 
Section~6.5}) of $\Cal T$, with $T$ intrepreted as the 
generic model in this topos.  For technical reasons, 
however, it is easier to work with the classifying topos of 
a slightly stronger theory, $\Cal T'$, obtained from $\Cal 
T$ by adding the following axiom for every term $t$ that 
contains the variable $x$ but is not just $x$.
\roster
\item[5] $t=x\implies\false$
\endroster
This additional axiom schema says that no (generalized) 
tree in a model of $\Cal T'$ can be a proper subtree of 
itself.  It is satisfied by the initial algebra (of finite trees), 
but not by, for example, the algebra of all finite and 
infinite binary trees, where the full binary tree (every 
node of which has two children) satisfies $x=[x,x]$.

Notice that, apart from making the proof easier, the 
addition of axiom schema \therosteritem5 to $\Cal T$ 
slightly improves the theorem.  The theorem remains 
true (with the same proof) if we replace $\Cal L$ by the 
stronger theory $\Cal L'$ where $T$ is assumed to satisfy 
\therosteritem5.  This strengthening of $\Cal L$ weakens 
the hypothesis of Theorem~4 and thus strengthens the 
theorem. 

As indicated above, we shall work with the classifying 
topos $\Cal E$ of $\Cal T'$.  In it, there is a generic (or 
universal) model $G$ of $\Cal T'$; this implies that, for 
any model $M$ of $\Cal T'$ in any Grothendieck topos 
$\Cal F$, there is a geometric morphism $\mu:\Cal F\to 
\Cal E$ whose inverse image functor $\mu^*$ sends $G$ 
to $M$.  (It implies more than this, but this, along with 
the explicit construction described below, will suffice for 
our purposes.)  Being a Grothendieck topos, $\Cal E$ gives 
an interpretation of higher order logic (as in \cite {4}) 
and local set theories (\cite {1}), with natural number 
object,  and intuitionistic Zermelo-Fraenkel set theory 
(\cite {5}), so by interpreting $T$ as $G$ we obtain a 
model of $\Cal L$ in the internal logic of $\Cal E$.  By the 
hypothesis of Theorem~4, it must be internally true in 
$\Cal E$ that there is a bijection from $P(G)$ to $Q(G)$.

The next part of the proof consists of studying $\Cal E$ 
and $G$ in sufficient detail to draw useful conclusions 
from this internal information.  We begin by describing 
$\Cal E$ explicitly as the topos of sheaves over a specific 
site.  As explained in \cite {9}, it is convenient to first 
build the classifying topos for the universal Horn axioms 
and then obtain $\Cal E$ as a sheaf subtopos.  For $\Cal 
T'$, axioms \therosteritem1, \therosteritem2, 
\therosteritem3, and \therosteritem5 are universal Horn 
sentences.  The classifying topos for this subtheory is, 
according to \cite{2}, the topos of presheaves on the dual 
of the category $\Cal A$ of finitely presented models of 
\therosteritem1, \therosteritem2, \therosteritem3, and 
\therosteritem5.  So we need to analyze such models 
$\langle A\mid E\rangle$, where $A$ is a finite set of 
generators and $E$ is a finite, consistent set of equations 
between terms built from the generators, 0, and $[-,-]$.

We show first that every such model is free, i.e., is 
isomorphic to $\langle B\rangle= \langle B\mid 
\emptyset\rangle$ for some finite set $B$.  To see this, 
we systematically simplify the given set $E$ of equations 
as follows.  (Technically, the simplification is an inductive 
process; at each step, we either decrease the cardinality 
of $A$ or we leave this cardinality unchanged but 
decrease the total length of all the equations in $E$.)  If 0 
occurs as one side of an equation in $E$, then the other 
side must be 0 or a member of $A$; it cannot be of the 
form $[t_1,t_2]$ because then the equation would be 
inconsistent by \therosteritem1.  If it is 0, then the 
equation $0=0$ can be deleted from $E$ as it is always 
true.  If it is a member $a$ of $A$, then we can delete $a$ 
from $A$, delete the equation $0=a$ (or $a=0$) from $E$, 
and replace all occurrences of $a$ in the rest of $E$ by 0.  
The result is a simpler presentation of the same algebra.  
So we may assume from now on that 0 does not occur as 
a side of an equation in $E$.  Suppose next that an 
element $a$ of $A$ occurs as a side of an equation in $E$.  
If the other side is also $a$, then the equation $a=a$ can 
simply be omitted.  Otherwise, the other side must be a 
term $t$ not involving $a$, for if it involved $a$ then the 
equation would contradict \therosteritem5.  So we can 
delete $a$ from $A$, delete $a=t$ (or $t=a$) from $E$, and 
replace all other occurrences of $a$ in $E$ by $t$.  Again, 
we have a simpler presentation of the same algebra.  So 
we may assume that each equation in $E$ has the form 
$[t_1,t_2]=[t_3,t_4]$; but such an equation can, thanks to 
\therosteritem2 and \therosteritem3, be replaced with 
the two equations $t_1=t_3$ and $t_2=t_4$, of lesser total 
length.  So again we get a simpler presentation of the 
same algebra.  Repeating these steps, we find that the 
process must terminate, for the size of $A$ cannot 
decrease infinitely often, and, after it stops decreasing, 
the total length of $E$ cannot decrease infinitely often.  
But the only way the process can stop is if $E$ has 
become empty.  This proves that $\langle A\mid 
E\rangle$ is isomorphic to $\langle 
B\mid\emptyset\rangle$ for some $B$ (a subset of $A$).

By virtue of this simplification, we may regard $\Cal A$ 
as consisting of only the free algebras $\langle A\rangle$  
on finite sets $A$ of generators.  We may also suppose 
that the only sets $A$ occurring are of the form 
$\{1,2,\dots,k\}$ for natural numbers $k$, since every 
$A$ is isomorphic to one of these.  We write $\a k$ for 
$\a{\{1,2,\dots,k\}}$.

The elements of $\a k$ are the variable-free terms of the 
language having the constant symbol 0, the binary 
operation $[-,-]$, and constant symbols for the generators 
$1,2,\dots,k$.  They can be identified with trees in which 
leaves may (but need not) be labeled with integers in the 
range from 1 to $k$.  So they are like patterns (defined in 
Section~1) except for the restriction on the possible labels 
and the fact that several leaves are allowed to have the 
same label.  We call them $k$-{\sl labeled trees\/}.  Note 
in particular that the members of $\a 0$ are simply the 
trees. 

A morphism in $\Cal A$ from $\a k$ to $\a l$ is, since 
$\a k$ is free, simply a map from $\{1,2,\dots,k\}$ into 
$\a l$, i.e., a $k$-tuple of $l$-labeled trees.  To compose 
this with some $\a l\to\a m$, i.e., with an $l$-tuple of 
$m$-labeled trees, take the $k$-tuple of $l$-labeled trees 
and replace, in each of its component trees, each leaf 
labeled $j$ with the $j$th $m$-labeled tree in the given 
$l$-tuple.  The identity morphism of $\a k$ is the 
$k$-tuple whose $i$th member is a single node labeled 
$i$.

The topos of set-valued functors on $\Cal A$ is the 
classifying topos for the universal Horn theory 
axiomtized by \therosteritem1, \therosteritem2, 
\therosteritem3, and \therosteritem5.  The universal 
model $U$ is the underlying set functor. For details about 
this, see for example \cite{2}.

To obtain the classifying topos for the full theory $\Cal 
T'$, we must pass to the subtopos of sheaves for the 
Grothendieck topology ``forcing'' the remaining axiom, 
\therosteritem4.  See \cite{9, 11} for more information 
about forcing topologies.  In the case at hand, the topology 
in question is that described in Part \therosteritem1 of 
the following lemma, whose other parts give useful 
alternative ways of viewing this topology.  In connection 
with Parts \therosteritem2 and \therosteritem3 of the 
lemma, recall that in the proof of Theorem~3 we 
introduced, for any polynomial $P(X)$, a set $S_n(P)$ of 
tagged tuples of patterns; we shall need this for the 
polynomials $X^k$.  In this special case of monomials, all 
tags are 1, so they can be omitted, and all the tuples are 
$k$-tuples.  So the elements of $S_n(X^k)$ can be taken to 
be simply $k$-tuples of patterns and thus, by suitable 
choice of labels, morphisms $\a k\to\a l$ for certain $l$'s.

\proclaim{Lemma  }
The following four Grothendieck topologies on the dual of 
$\Cal A$ coincide.
\roster
\item The smallest topology for which $\a1$ is covered 
by the set of two morphisms $\a1\to\a0:1\mapsto0$ and 
$\a1\to\a2:1\mapsto[1,2]$
\item The smallest topology in which, for each $n$, each 
$\a k$ is covered by the set $S_n(X^k)$.
\item The topology where the covering sieves of any $\a 
k$ are those sieves that include $S_n(X^k)$ for some $n$.
\item The topology where the covering sieves of any $\a 
k$ are those sieves that include a finite family of maps 
$\a k\to\a {l_i}$ such that every map $\a k\to\a0$ 
factors through a map from the finite family.
\endroster
(It is part of the assertion of the lemma that the 
collections of sieves described in \therosteritem3 and 
\therosteritem4 are topologies.)
\endproclaim

\demo{Proof}
We first show that the family of sieves in described in 
\therosteritem4 is a Grothendieck topology.  It clearly 
contains the maximal sieve on any object (use the family 
consisting of just the identity map).  If it contains a sieve 
$R$ on $\a k$, witnessed by a finite family of $\a k\to\a 
{l_i}$ as in \therosteritem4, and if $f:\a k\to\a m$ is any 
morphism, then those pushouts of the $\a k\to\a 
{l_i}$ along $f$ that exist in $\Cal A$ witness that 
\therosteritem4 also contains the sieve of morphisms out 
of $\a m$ whose composites with $f$ are in $R$, i.e., the 
pullback of $R$ along $f$ in the sense of the dual 
category.  (We use here that, for any pair of maps in 
$\Cal A$ with the same domain, if they can be completed 
to a commutative square then they have a pushout.  This 
is a general property of categories of models of universal 
Horn theories.)  Finally, the alleged topology is closed 
under composition, because if a family of maps $\a 
k\to\a {l_i}$ and, for each $i$, another family of maps $\a 
{l_i}\to\a {m_{ij}}$ satisfy the requirements in 
\therosteritem4, then so does the family of all $\a k\to\a 
{m_{ij}}$.  Thus, \therosteritem4 is a topology.

This topology contains the sieve generated by the two 
maps in \therosteritem1.  Indeed, these two maps 
themselves serve as the finite family required in 
\therosteritem4, for any morphism $\a1\to\a0$ must 
send 1 to either the empty tree 0 or a non-empty tree 
$[t_1,t_2]$, and in the former case it factors through the 
specified map $\a1\to\a0$ while in the latter case it 
factors through the specified map $\a1\to\a2$ (via the 
map sending 1 and 2 to $t_1$ and $t_2$).  Therefore, the 
topology \therosteritem4 includes the topology 
\therosteritem1.

We show next that \therosteritem1 includes 
\therosteritem2; of course it suffices to show that the 
generating covers $S_n(X^k)$ in \therosteritem2 are also 
covers with respect to \therosteritem1.  For this purpose, 
note that $S_n(X^k)$ can be obtained from the identity 
map of $\a k$ (a $k$-tuple of distinctly labeled, 
one-element trees) by repeated use of the development 
process described in the proof of Theorem~3.  That is, a 
tree $r$ with a particular labeled node is replaced by two 
trees, $r'$ where that node has been removed and $r''$ 
where that node has been given two distinctly labeled 
children and has lost its own label.  Notice that, in each of 
the $k$-tuples arising in this development process, no 
labels are repeated.  We show that applying one 
development step to a cover in the topology 
\therosteritem1, with no repeated labels in its 
$k$-tuples, yields again a cover.  This will clearly imply 
that $S_n(X^k)$, obtained by repeated development of a 
trivial cover, is itself a cover, as desired.  So consider one 
development step applied to such a cover.  Suppose it 
replaces a node labeled $i$ in some (unique)  component 
of a $k$-tuple, $\a k\to\a l$, in the given covering.  Then 
the pushouts in $\Cal A$ of the two maps in 
\therosteritem1 along $\a1\to\a l:1\mapsto i$ cover $\a 
l$ (since a topology on the dual of $\Cal A$ is closed 
under pullbacks in this dual).  So in the given cover of 
$\a k$ we may replace the map $\a k\to\a l$ by its 
composites with these two pushouts, and we still have a 
cover of $\a k$.  But this replacement is precisely the 
development step at label $i$.  This completes the proof 
that topology \therosteritem1 includes \therosteritem2.

It is trivial that all the sieves described in 
\therosteritem3 are in topology \therosteritem2.  That 
\therosteritem3 is a topology will follow once we show 
that it includes \therosteritem4, for then all four items 
listed in the lemma are equal.

So, to complete the proof, we consider an arbitrary sieve 
$R$ containing a finite family of maps $f_i:\a k\to\a {l_i}$ 
as in \therosteritem4, and we show that this sieve is in 
\therosteritem3.  Fix an integer $n$ greater than the 
depths of all the labeled trees occurring in the $k$-tuples 
$f_i$.  We would like to develop all these $k$-tuples to 
depth $n$ and show that the resulting $k$-tuples are all 
in $R$ and include all of $S_n(X^k)$.  Some caution is 
needed, however, since the same label may occur several 
times in an $f_i$, and development has not even been 
defined for such an $f_i$.  We begin by extending the 
notion of development to the case of repeated labels.  If 
label $z$ occurs several times in a $k$-tuple $\vec r$ 
(either in the same component tree or in different 
components) and if all its occurrences are at depths $\leq 
n$, then we develop $\vec r$ at $z$ by replacing it by 
two $k$-tuples, $\overarrow{r'}$ where all nodes labeled 
$z$ have been deleted, and $\overarrow{r''}$ where each 
node labeled $z$ has been (unlabeled and) given two 
children, the left one being labeled $z_l$ and the right 
one $z_r$, where these are two labels not yet occurring in 
$\vec r$.  In other words, we carry out development in 
the previous sense at all $z$-labeled nodes in parallel, 
treating them all identically.  As in the proof that 
\therosteritem1 includes \therosteritem2, it is easy to 
see that this definition makes $\overarrow{r'}$ and 
$\overarrow{r''}$, regarded as morphisms out of $\a k$, 
the composites of $\vec r$ with two other morphisms; in 
particular, $\overarrow{r'}$ and $\overarrow{r''}$ belong 
to the sieve $R$ if $\vec r$ does.  By repeated 
development, the given family of maps $f_i$ becomes a 
new family of maps $f'_j$ with the property that each 
tree in it has depth at most $n+1$, and every label that 
occurs in any $f'_j$ has at least one occurrence at depth 
$n+1$ in $f_j$ (for if all occurrences were at depth $\leq 
n$ then we could develop further at that label).  After 
development is finished, any $f'_j$ without repeated 
labels is a member of $S_n(X^k)$.

We shall show that all members of $S_n(X^k)$ must be 
among the $f'_j$'s; as pointed out above, this will suffice 
to complete the proof, since each $f'_j$ is in $R$.  
Suppose, toward a contradiction, that $\vec p\in 
S_n(X^k)$ is not among the $f'_j$'s.  Form an instance 
$\vec t$ of $\vec p$ by a substitution that replaces the 
labels in $\vec p$ with distinct trees of all of the same 
depth $d>n$.  This $\vec t$ is an instance of at least one 
of the original $f_i$'s, since all $k$-tuples of trees are 
such instances.  By considering the development process 
one step at a time, one easily sees that $\vec t$ is an 
instance of some $k$-tuple at each stage of the 
development; in particular it is an instance of some $f'_j$ 
at the final stage.  That $f'_j$ cannot have distinct labels, 
for then it would be in $S_n(X^k)$, whereas the unique 
element of $S_n(X^k)$ having $\vec t$ as an instance is 
$\vec p$ which is not an $f'_j$.  So $\vec t$ is an instance 
of an $f'_j$ in which some label, say $z$, occurs at least 
twice.  Let $s$ be the tree substituted for $z$ in 
instantiating $f'_j$ to $\vec t$.  Then $s$ occurs at least 
twice as a subtree in components of $\vec t$, namely 
wherever $z$ occurred as a label in $f'_j$, and at least 
one of these subtrees has its root at depth $n+1$ in $\vec 
t$.  So when $\vec t$ was obtained by instantiating $\vec 
p$, some label must have been replaced by $s$; thus, by 
our choice of that substitution, $s$ has depth $d>n$.  Our 
choice of substitution also ensures that $s$ cannot have a 
second occurrence with its root at depth $n+1$ in $\vec t$, for 
at depth $n+1$ we replaced all nodes in $\vec p$ with 
different trees.  Nor can a second occurrence of $s$ in 
$\vec t$ have its root at depth greater than $n+1$, for 
then that occurrence would extend to depth greater than 
$n+d$, which is the maximum depth occurring in $\vec t$.  
So $s$ must have a second occurrence in $\vec t$ with its 
root $x$ at depth $e\leq n$.  In $\vec p$, there is a node 
at the same position that $x$ has in $\vec t$ (as the 
substitution leading from $\vec p$ to $\vec t$ affects 
only depths $>n$), and we call it $x$ also.  If $x$ has any 
labeled descendants  in $\vec p$, then in $\vec t$ those 
descendants, at level $n+1$, were replaced by trees of 
depth $d$.  So the subtree with root $x$ in $\vec t$ has 
depth $n-e+1+d>d$.  If, on the other hand, $x$ has no 
labeled descendants in $\vec p$, then the subtree with 
root $x$ extends to depth at most $n-1<d$ in $\vec p$ (as 
$\vec p\in S_n(X^k)$, so any node at depth $n$ must 
have two labeled children), and nothing changes in this 
subtree when we pass from $\vec p$ to $\vec t$.  But the 
subtree of $\vec t$ with root $x$ is $s$, of depth $d$, so 
both cases are contradictory.  This contradiction 
completes the proof that the sieve $R$ includes 
$S_n(X^k)$ and is 
therefore a covering in \therosteritem3.
\qed\enddemo

Let $\Cal E$ be the topos of sheaves on the dual of $\Cal 
A$ with respect to the topology described in this lemma.  
We write $i$ for the canonical geometric morphism from 
$\Cal E$ to the presheaf topos $\Cal S^{\Cal A}$, so $i_*$ 
is the inclusion functor and $i^*$ is the associated sheaf 
functor.  It follows from \cite{9} that $\Cal E$ is the 
classifying topos for models of $\Cal T'$, the generic 
model being $G=i^*(U)$.

The hypothesis of the theorem implies that the statement 
``there exists a pair of inverse bijections $f:P(G)\to Q(G)$ 
and $g:Q(G)\to P(G)$'' is internally true in $\Cal E$.  This 
implies, by virtue of the internal meaning of ``there 
exists,'' that there is an epimorphism $C\to1$ in $\Cal E$ 
such that in the slice topos $\Cal E/C$ there is a pair of 
actual inverse isomorphisms between $C^*(P(G))$ and 
$C^*(Q(G))$, where $C^*$ means the inverse image along 
the canonical geometric morphism $\Cal E/C\to\Cal E$.  

We show that $C$ has a global element $1\to C$.  The fact 
that $C\to 1$ is an epimorphism in the sheaf topos $\Cal 
E$ means that, in the presheaf topos, it has dense range, 
i.e., that every object of the site is covered by arrows 
from objects $A$ where $C(A)$ is inhabited.  (Here 
``from'' refers to the site, the dual of $\Cal A$; in terms of 
maps in $\Cal A$ we would have ``to''  instead.)  But 
$\a0$ is covered only by its maximal sieve, as is clear by 
description \therosteritem4 in the lemma.  So $C(\a0)$ is 
inhabited, say by $z$.  But $\a0$ is initial in $\Cal A$, so 
$C$ has a global section (in $\Cal S^{\Cal A}$ and hence 
also in $\Cal E$) whose value at any object $A$ is the 
image of $z$ under $C$ of the unique map $\a0\to A$.

By taking inverse images along this section, we find that 
already in $\Cal E$ (without having to pass to a slice 
topos) we have a pair of inverse isomorphisms  $f:P(G)\to 
Q(G)$ and $g:Q(G)\to P(G)$.

Now consider, in the topos $\Cal S$ of sets, the set $T$ of 
(finite, binary) trees.  It is a model of $\Cal T'$ in an 
obvious way (in fact the initial model of $\Cal T'$).  So it 
is $\mu^*(G)$ for some geometric morphism $\mu^*:\Cal 
S\to\Cal E$ (a point of $\Cal E$).  Since inverse images 
along geometric morphisms preserve finite products and 
coproducts, $\mu^*(f)$ and $\mu^*(g)$ are inverse 
bijections between the sets $P(T)$ and $Q(T)$.  To 
complete the proof of the theorem, it suffices to show 
that they are very explicit, for then the polynomials 
$P(X)$ and $Q(X)$ are combinatorially equivalent and, by 
Theorem~3, algebraically equivalent, as desired.

In fact, it suffices to prove somewhat less.  In the 
definition of ``very explicit function'' in Section~3, we 
required that no label be repeated in any $\vec q_i$ (as 
part of the notion of pattern) and that $\vec q_i$ contain 
exactly the same labels as $\vec p_i$.  But we mentioned, 
just before Theorem~3, that one could relax these 
requirements, by allowing repeated labels in $\vec q_i$ 
and by allowing labels to occur in $\vec p_i$ without 
occurring in $\vec q_i$, without affecting the very 
explicit bijections.  This is because any function that is 
very explicit in the liberalized sense but not in the 
original sense cannot be a bijection.  So in our present 
situation, it suffices to show that $\mu^*(f)$ and 
$\mu^*(g)$ are very explicit in the liberalized sense, since 
we already know that they are bijections.

Of course it suffices to treat $\mu^*(f)$, as the situation is 
symmetric between the two bijections.  In fact, it suffices 
to treat the restriction of $\mu^*(f)$ to one of the 
summands $T^k$ in $P(T)$, for if each of these 
restrictions is very explicit then so is $\mu^*(f)$ itself.  
Note that such a restriction of $\mu^*(f)$ is $\mu^*$ of a 
restriction of $f$ to one of the summands $G^k$ of $P(G)$ 
in $\Cal E$.

So, changing notation slightly, we have a morphism 
$f:G^k\to Q(G)$ in $\Cal E$, and we wish to show that 
$\mu^*(f):T^k\to Q(T)$ is very explicit in the liberalized 
sense.

Recall that $G$ is obtained by applying the associated 
sheaf functor $i^*$ to the underlying set functor $U$ in 
$\Cal S^{\Cal A}$.  As $i^*$ preserves sums and products, 
we can write $f:i^*(U^k)\to i^*(Q(U))$, so $f$ corresponds, 
under the adjunction $i^*\dashv i_*$, to a map of 
presheaves $U^k\to i_*i^*(Q(U))$. Since $U^k$ is a 
representable presheaf, represented by $\a k$, such a 
map corresponds via Yoneda's Lemma to an element of 
$(i_*i^*(Q(U)))(\a k)$.  Here $i_*i^*(Q(U))$ is the associated 
sheaf of $Q(U)$, regarded as a presheaf.  Inspecting the 
usual construction of associated sheaves in terms of 
``patched together'' sections, we find that any element of 
$(i_*i^*(Q(U)))(\a k)$ can be described as follows.  There 
is a cover of $\a k$, say by maps $\a k\to \a{l_i}$, and to 
each of these maps is assigned an element $q_i\in 
(Q(U))(\a {l_i})=Q(U(\a {l_i}))$ in a coherent manner.  
(Coherence means that if two of the maps $\a k\to 
\a{l_i}$ in the cover can be completed to a commutative 
square by maps $\a{l_i}\to A$, then the resulting images 
in $Q(U(A))$ of the $q_i$'s are equal.)
By the lemma, we may take the covering to be $S_n(X^k)$ 
for some $n$; then the maps $\a k\to \a{l_i}$ in the 
covering are $k$-tuples of patterns (since no label is 
re-used in a map in $S_n(X^k)$) such that every 
$k$-tuple of trees is an instance of exactly one of them.  
If we call these patterns $p_i$, then they and the 
corresponding $q_i$ are exactly as required in the 
liberalized definition of a very explicit function.  It 
remains to check that, if we start with $f:G^k\to Q(G)$, 
transform it to $U^k\to i_*i^*(Q(U))$ and then to an 
element of $(i_*i^*(Q(U)))(\a k)$, and finally represent 
that element on a covering to obtain a very explicit 
function, as just described, then that function agrees with 
$\mu^*(f)$.  This verification is routine, tedious, and 
therefore omitted.
\qed\enddemo

Lawvere has pointed out that the proof of Theorem~4 emphasizes the
distinction between the 2-categorical sort of universality that
defines classifying topoi and the 1-categorical sort that defines free
algebras.  For the proof shows that the (2-categorically universal)
generic model $G$ of $\Cal T'$ satisfies $P(G)\cong Q(G)$ only when
$P(X)$ and $Q(X)$ are algebraically equivalent; in particular, $G^2$
and $G$ are not isomorphic in the classifying topos.  In contrast, as
we remarked before the proof of Theorem~4, the (1-categorically
universal) free model $T$ of $\Cal T'$ in any topos, consisting of
finite trees, does satisfy $T^2\cong T+T\cong T$.

\head
Acknowledgment
\endhead

I thank Bill Lawvere for helpful conversations and e-mail messages
about the subject of this paper.


\Refs

\ref\no1
\by J. L. Bell
\book Toposes and Local Set Theories
\bookinfo Oxford Logic Guides 14
\publ Clarendon Press
\publaddr Oxford
\yr1988
\endref

\ref\no2
\by A. Blass and A. Scedrov
\paper Classifying topoi and finite forcing
\jour J. Pure Appl. Alg.
\vol 28
\yr 1983
\pages 111--140
\endref



\ref\no3
\by S. Burris and H. P. Sankappanavar
\book A Course in Universal Algebra
\bookinfo Graduate Texts in Mathematics 78
\publ Springer-Verlag
\yr1981
\endref

\ref\no4
\by M. P. Fourman
\paper The logic of topoi
\inbook Handbook of Mathematical Logic
\ed J. Barwise
\publ North-Holland
\yr1977
\pages1053--1090
\endref

\ref\no5
\by M. P. Fourman
\paper Sheaf models of set theory
\jour J. Pure Appl. Alg.
\vol 19
\yr 1980
\pages 91--101
\endref

\ref\no6
\by A. Garsia and S. Milne
\paper Method for constructing bijections for classical partition 
identities
\jour Proc. Nat. Acad. Sci. U. S. A.
\vol 78
\yr 1981
\pages 2026--2028
\endref

\ref\no7
\by P. T. Johnstone
\book Topos Theory
\bookinfo London Mathematical Society Monographs 10
\publ Academic Press
\yr1977
\endref



\ref\no8
\by F. W. Lawvere
\paper Some thoughts on the future of category theory
\inbook Category Theory, Proceedings, Como 1990
\bookinfo Lecture Notes in Mathematics 1488
\eds A. Carboni, M. C. Pedicchio, and G. Rosolini
\publ Springer-Verlag
\yr1991
\pages1--13
\endref

\ref\no9
\by A. Scedrov
\book Forcing and Classifying Topoi
\bookinfo Memoirs Amer. Math. Soc. 295 
\yr1984
\endref

\ref\no 10
\by S. Schanuel
\paper Negative sets have Euler characteristic and dimension
\inbook Category Theory, Proceedings, Como 1990
\bookinfo Lecture Notes in Mathematics 1488
\eds A. Carboni, M. C. Pedicchio, and G. Rosolini
\publ Springer-Verlag
\yr1991
\pages379--385
\endref



\ref\no11
\by M. Tierney
\paper Forcing topologies and classifying topoi
\inbook Algebra, Topology and Category Theory: a collection of 
papers in honor of Samuel Eilenberg
\ed A. Heller and M. Tierney
\publ Academic Press
\yr1976
\pages211--219
\endref

\endRefs
\enddocument